\documentclass[11pt,oneside]{article}
\usepackage{amsmath,amsthm,amsfonts,amssymb,MnSymbol,mathrsfs,latexsym}

\usepackage[english]{babel}
\usepackage{graphicx}

\usepackage{a4wide}

\usepackage[applemac]{inputenc}
\usepackage[all]{xy}

\numberwithin{equation}{section}

\usepackage{bm}

  \newcommand{\const}{\rm const}
  \newcommand{\Var}{\rm Var}
  
  \newcommand{\Law}{\rm Law}

\theoremstyle{plain}
\newtheorem{theorem}{Theorem}[section]
\theoremstyle{theorem}
\newtheorem{corollary}{Corollary}[section]
\newtheorem{lemma}[theorem]{Lemma}
\newtheorem{proposition}[theorem]{Proposition}

\newtheorem{remark}{Remark}[section]
\newtheorem{example}{Example}[section]

\renewenvironment{proof}{{\bf{Proof.}}}{\hfill $\Box$ \\}

\usepackage{lipsum}
\usepackage{fancyhdr}

\title{\large \textbf{Gaussian and non - Gaussian distributed \\[1mm] random  analytical and entire functions}}

\footnotesize\date{}

\author{\normalsize \textbf{Maria Rosaria FORMICA ${}^{1}$},   \normalsize\textbf{Eugeny OSTROVSKY
${}^2$} and \normalsize\textbf{Leonid SIROTA ${}^3$}}

\begin{document}

\maketitle

\begin{center}
{\footnotesize ${}^{1}$ Universit\`{a} degli Studi di Napoli \lq\lq Parthenope\rq\rq, via Generale Parisi 13,\\
Palazzo Pacanowsky, 80132,
Napoli, Italy.} \\

\vspace{2mm}

{\footnotesize e-mail: mara.formica@uniparthenope.it} \\

\vspace{4mm}

{\footnotesize ${}^{2,\, 3}$  Bar-Ilan University, Department of Mathematics and Statistics, \\
52900, Ramat Gan, Israel.} \\

\vspace{2mm}

{\footnotesize e-mail: eugostrovsky@list.ru}\\

\vspace{2mm}

{\footnotesize e-mail: sirota3@bezeqint.net} \\

\end{center}

 \vspace{3mm}

\begin{abstract}

We investigate the complex Gaussian as well as non-Gaussian
distributed random analytical and entire functions (complex entire
random field) and calculate their domain of definiteness (radius of
convergence) as well as some important characteristics: order and
type.

As a consequence we deduce that all the mentioned characteristics,
under very natural conditions, are deterministic (non-random) with
probability one and we calculate them. Moreover we exhibit some
examples to show the exactness of the obtained results.

\end{abstract}

\vspace{4mm}

 {\it \footnotesize Keywords:}
{ \footnotesize Gaussian and non-Gaussian distributed random
variables, entire functions, Taylor (power) series, radius of
convergence, tail function, centered variables and random functions,
Lemma of Borel-Cantelli, complex variables and functions, order and
type of an entire function,  distribution of zeros.}

\vspace{2mm}

\noindent {\it \footnotesize 2010 Mathematics Subject
Classification}:
{\footnotesize
30D20,   
30D10,    
30B20    
}

\vspace{4mm}

\section{Introduction.}

\vspace{2mm}

Define the following random entire Gaussian distributed  centered
complex function of the form (power, Taylor series)
\begin{equation} \label{defin fun series}
f(z) = f(z,\omega) = f[\{\xi_k\}](z) \stackrel{def}{=} \sum_{k=0}^{\infty} \xi_k \ z^k.
\end{equation}
Here $ \ z \ $ is a {\it complex} variable, $ \ \{\xi_k\}, \ k =
0,1,2,\ldots \ $, are centered (mean zero) {\it Gaussian
distributed} complex random variables (r.v.), defined on some
sufficiently rich probability space $ \ \Omega = (\{\omega\}, B,
{\bf P}): \ $
$$
\xi_k = \eta_k + i \ \zeta_k, \ \ i^2 = - 1,
$$
not necessarily independent and, in general case, the variables$ \
\eta_k, \ \zeta_k, \ $ are also dependent. Here the random variables
$ \ \eta_k, \zeta_k \ $ are real valued. \par

In Section \ref{Sec Gaussian case} we analize the Gaussian case
while in Section \ref{section
 non-gaussian case} we consider the non - Gaussian case.

 \ Denote also

\begin{equation} \label{complex case 1}
\beta_k^2 = \Var \ \eta_k^2 = {\bf E} \eta_k^2, \ \ \gamma_k^2 =
\Var \ \zeta_k^2 = {\bf E} \zeta_k^2,
\end{equation}

\begin{equation} \label{complex case 2}
\sigma_k^2 = \beta_k^2 + \gamma_k^2.
\end{equation}

\vspace{4mm}

 In this short report we study the following characteristics of the introduced function: radius of convergence,  order
and type, etc. \vspace{4mm}

 \ Some properties of these random entire functions are studied in many works, see e.g.
 \cite{Cheng,Gao,Hough,Nazarov,Peres,Sodin1,Sodin2}.

 Throughout this paper the letters $ \ C, C_j(\cdot) \ $ will denote, as ordinary, various non-essential positive finite
constants which may differ from one formula to the next even within
a single string  of estimates and which does not depend on the
 variables.  We make no attempts to obtain the best values
for these constants.

\vspace{4mm}

\section{Main results. Gaussian entire random function.}\label{Sec
Gaussian case}

\begin{center}

\vspace{4mm}

 \ {\sc  Radius of convergence. Gaussian case.} \par

\vspace{4mm}

\end{center}

\begin{theorem}\label{Theorem radius convergence}
 The radius $ \ R \ $  of convergence of the Taylor series
\eqref{defin fun series} is non-random
 a.e. and  may be calculated by the formula

\begin{equation} \label{radii}
R = \frac{1}{\overline{\lim}_{n \to \infty}   \sqrt[n]{\sigma_n}}.
\end{equation}

\end{theorem}

The proof of Theorem \ref{Theorem radius convergence} is a
consequence of the following auxiliary important lemma.

%
%
%
%

 \begin{lemma}\label{lemma limsup}
  Let $ \ \{\xi_n \}, \ n = 1,2,\ldots \ $, be a sequence of centered Gaussian distributed
  random variables. We state that, with probability one,

\begin{equation} \label{over lim}
\overline{\lim}_{n \to \infty}   \sqrt[n]{|\xi_n|}  =
\overline{\lim}_{n \to \infty}   \sqrt[n]{\sigma_n}.
\end{equation}

\end{lemma}

\vspace{3mm}


 \begin{proof}

\vspace{2mm}
 {\it Upper estimate.}

 \vspace{2mm}

 We follow the arguments in \cite[Lemma 7]{Cheng}. Suppose, without loss of generality, that
$$
\overline{\lim}_{n \to \infty}  \sqrt[n]{\sigma_n} = 1.
$$

 \ Let us consider at first the case when $ \  \xi_n \ $ are numerical valued variables. Let $ \ \epsilon > 0 \ $  be given. There exists an integer
 number $ \ n_0 = n_0(\epsilon) \ $ for which

$$
\sigma_n < (1 + \epsilon)^n,  \ \ n \ge n_0.
$$

 \ The probability $ \ {\bf P_n}(\epsilon) := {\bf P}  \left(  \sqrt[n]{ |\xi_n|} > (1 + \epsilon)^2   \right) \ $
allows the following estimate

\begin{equation} \label{prob estim}
{\bf P_n}(\epsilon) \le 2 \exp \left( \ - 0.5 (1 + \epsilon)^n \   \right).
\end{equation}

 \ Since
$$
\forall \epsilon > 0 \ \Rightarrow \sum_{n=1}^{\infty} {\bf P_n}(\epsilon) < \infty,
$$
we deduce, by virtue of the lemma of Borel - Cantelli, that

$$
\overline{\lim}_{n \to \infty} \sqrt[n]{|\xi_n|} \le 1  \ a.e.
$$

\vspace{5mm}

 \ {\it Lower estimate.}

 \vspace{2mm}

  \  One can suppose, without loss of generality, $ \ \sigma_n = 1, \  n \ge 1. $
  Let again $ \ \epsilon = \const \in (0,1). \ $ Consider the probability

 $$
 {\bf P}^{(n)}(\epsilon) \stackrel{def}{=} {\bf P} \left( \sqrt[n]{|\xi_n|} \le 1 - \epsilon  \right).
 $$
 \ We have

$$
{\bf} P^{(n)}(\epsilon) = {\bf P} \left( |\xi|_n \le (1 - \epsilon)^n   \right) \le C \ (1 - \epsilon)^n,
$$
so that, as before,

$$
\forall \epsilon \in (0,1) \ \Rightarrow  \sum_{n=1}^{\infty} {\bf} P^{(n)}(\epsilon)  < \infty.
$$
 \ Therefore,

$$
\overline{\lim}_{n \to \infty} \sqrt[n]{|\xi_n|} \ge 1 \ a.e.
$$

 Let us return to the general case, i.e. the complex case \eqref{complex case 1}, \eqref{complex case 2}. As above,
$ \ \xi_n = \eta_n + i \ \zeta_n, \ i^2 = - 1, $ so that the value $ \  ||\xi_n||^2  \ $ may be represented in the form

$$
||\xi_n||^2 = \beta_n^2 \ \eta_n^2 + \gamma_n^2 \ \zeta_n^2,
$$
where $ \ \Law \{\eta_n\} =  \Law \{\zeta_n\} = N(0,1), \ $ and one
can suppose, without loss of generality,

$$
\beta_n^2 + \gamma_n^2 = \sigma_n^2 = 1.
$$

 \ Denote, as usually, by $ \ ||\tau||_p \ $  the ordinary  Lebesgue - Riesz $ \ L_p \ $ norm of the random variable $ \ \tau: \ $
$$
||\tau||_p \stackrel{def}{=} \left[ \ {\bf E} |\tau|^p \
\right]^{1/p}, \ \ p \ge 1.
$$
 \ It is well known that if $ \ \Law(\tau) = N(0,1), \ $ then

$$
||\tau||_p \asymp C_1 p^{1/2}, \ \ p \ge 1,
$$
therefore
$$
||\tau^2||_p \asymp C_2 \ p, \ p \ge 1.
$$

 \ We apply the triangle inequality  for the   Lebesgue - Riesz norm

 $$
 ||\xi_n||^2_p \le C_3 \ p,
 $$
or equally

$$
||\xi_n||_p \le C_4 \  \sqrt{p}, \ \ p \ge 1.
$$

 \ Hence, the r.v. - s  $ \ \{\tau_n\}, \ n = 1,2,\ldots \ $  are uniformly subgaussian:

$$
\sup_n {\bf P} (|\xi_n| > u) \le 2 \exp(- C_5 u^2), \ \ u \ge 0.
$$

 \ It follows immediately as before, from the last relation, that with probability one

$$
\overline{\lim}_{n \to \infty} \sqrt[n]{|\xi_n|} \le 1.
$$
 \ The inverse  relation

$$
\overline{\lim}_{n \to \infty} \sqrt[n]{|\xi_n|} \ge 1
$$
may be proved quite analogously to that in the one - dimensional
case.
\par

\vspace{5mm}

 \  Lemma \ref{lemma limsup} and consequently Theorem \ref{Theorem radius convergence}, are proved. \par


\end{proof}

 \ As a slight consequence, if

\begin{equation} \label{condit R}
 \lim_{n \to \infty} \sqrt[n]{\sigma_n} = 0,
\end{equation}

then $ \ R = \infty.\ $ \par

\vspace{4mm}

\begin{center}

 \ {\sc Order and type for the Gaussian entire functions.}

\end{center}

\vspace{4mm}

 \ Suppose in this subsection that the condition \eqref{condit R} is satisfied; then the  random function $ \ f = f(z) \ $ from
the definition of the Taylor (power) series \eqref{defin fun series}
is entire on the whole complex plane with probability one: $ \ R
=\infty. \ $ \par

 \ Recall the following classical definitions from the theory of an analytical complex  functions, see e.g. \cite[chapter 1]{Levin}.
 Let the complex valued entire function $ \ g = g(z) \ $ be defined (and be analytical)  on the whole complex plane. Denote

$$
M_g(r) := \sup_{z\,:\, |z| = r} |g(z)| = \sup_{z: |z| \le r} |g(z)|.
$$
 \ The value
$$
\rho[g] \stackrel{def}{=} \overline{\lim}_{r \to \infty}  \left\{ \frac{\ln \ \ln M_g(r) }{\ln r}   \right\}
$$
is said to be {\it Order} of the function $ \ g. \ $ The {\it Type} $  \beta[g] \ $ of these function $ \ g \ $ is defined
by an equality
$$
\beta[g] \stackrel{def}{=} \overline{\lim}_{r \to \infty} \left\{ \ \frac{\ln M_g(r)}{r^{\rho[g]}} \ \right\}.
$$

\vspace{3mm}

 \ Let us return to the source random complex function $ \ f = f(z) \ $ in \eqref{defin fun series}. \par

\vspace{5mm}

\begin{theorem}\label{theorem order-type non random}
 Both the characteristics for the random Gaussian function $ \ f = f(z) $, order and type, are a.e.
 non-random and may be calculated correspondingly by the relations

\vspace{4mm}

\begin{equation} \label{rho f}
\rho[f] = \overline{\lim}_{n \to \infty} \left\{ \ \frac{n \ \ln n}{|\ln \sigma_n \ |} \  \right\},
\end{equation}

\vspace{4mm}

\begin{equation} \label{beta f}
\beta[f] = \overline{\lim}_{n \to \infty} \left\{ \ n^{1/\rho[f]} \ \sqrt[n]{ \sigma_n} \  \right\}.
\end{equation}

\end{theorem}

\vspace{4mm}

 \begin{proof}  The equality \eqref{beta f} and in particular, the non-randomness of the value $ \ \beta[f] \
 $,
can be ground quite alike ones in the Theorem \ref{Theorem radius
convergence} and Lemma \ref{lemma limsup}.
\par

 \ Further, it is known, see e.g. \cite[chapter 1, pp. 5-15]{Levin}, that for the entire function of the form

$$
g(z) = \sum_{k=0}^{\infty} c_k  z^k
$$
there holds

\begin{equation} \label{rho g}
\rho[g] = \overline{\lim}_{n \to \infty} \left\{ \ \frac{n \ \ln n}{|\ln |c_n| \ |} \  \right\},
\end{equation}

\vspace{4mm}

\begin{equation} \label{beta g}
\beta[g] = \overline{\lim}_{n \to \infty} \left\{ \ n^{1/\rho[f]} \ \sqrt[n]{| c_n|} \  \right\}.
\end{equation}

 \ Therefore

\begin{equation} \label{rho rand f}
\rho[f] = \overline{\lim}_{n \to \infty} \left\{ \ \frac{n \ \ln n}{|\ln |\xi_n| \ |} \  \right\},
\end{equation}

\vspace{4mm}

\begin{equation} \label{beta rand f}
\beta[f] = \overline{\lim}_{n \to \infty} \left\{ \ n^{1/\rho[f]} \ \sqrt[n]{| \xi_n|} \  \right\}.
\end{equation}

 \ It is no hard to ground alike to the proof of Lemma \ref{lemma limsup} that
both the last expressions in the right - hand sides of equations
\eqref{rho rand f} and \eqref{beta rand f} are actually non-random
and coincide, correspondingly, with the ones in Theorem \ref{Theorem
radius convergence}. \par

\end{proof}
 \ We will ground in the next section a more general case, where the r.v. $ \ \{\xi_n\} \ $ are not necessarily
  Gaussian. \par

\vspace{5mm}

 \section{Non-Gaussian case.} \label{section
 non-gaussian case}

\vspace{2mm}

 Let us consider now the series of the form
\begin{equation} \label{defin non Gauss  series}
h(z) = h(z,\omega) = h[\{\eta_k\}](z) \stackrel{def}{=} \sum_{k=0}^{\infty} \eta_k \ z^k.
\end{equation}

 \ Here $ \ \{\eta_k\}, \ k = 0,1,2,\ldots  \ $,  are random variables, not necessarily independent, Gaussian or
centered.\par

\vspace{4mm}

 \ Introduce the following tail functions

\begin{equation} \label{tail functions}
T_k(u) := {\bf P}(|\eta_k| > u), \ \ \ S_k(u) := {\bf P}(|\eta_k|
\le u) = 1 - T_k(u), \ \ u \ge 0.
\end{equation}

\vspace{4mm}

\begin{theorem}\label{theorem non-Gaussian case}
 (Upper bound). Suppose that there exists a positive sequence of
{\it deterministic} numbers $ \ \sigma_k\ $ such that the {\it
uniform} tail function

\begin{equation} \label{tail non Gauss}
T(u) \stackrel{def}{=} \sup_{k} {\bf P} (|\eta_k|/\sigma_k > u), \ \
u \ge 1
\end{equation}
satisfies the following condition

\begin{equation} \label{T condition}
\forall Q > 1 \ \ \Rightarrow \ \ \sum_{k=1}^{\infty} T(Q^k) <
\infty.
\end{equation}
Then, with probability one,

\begin{equation} \label{Non Gauss radii}
\overline{\lim}_{n \to \infty} \sqrt[n]{|\eta_n|} \le  \overline{\lim}_{n \to \infty} \sqrt[n]{\sigma_n}.
\end{equation}

\vspace{4mm}

More generally assume that, for all the values $ \ Q > 1 \ $,

\begin{equation} \label{more general}
\sum_{k=1}^{\infty} {\bf P} (|\eta_k|/\sigma_k > Q^k) < \infty.
\end{equation}
 \ Then the relation  (\ref{Non Gauss radii}) again holds true by virtue of the well-known Lemma of Borel - Cantelli. \par
\end{theorem}
\vspace{4mm}

\begin{remark}
 {\rm The variables $ \ \sigma_k \ $ in \eqref{tail non Gauss} may
be chosen, for instance, as follows. Let $ \ ||\cdot|| \ $ be some
rearrangement invariant norm defined on spaces of random variables,
for example, Lebesgue-Riesz, Orlicz or Grand Lebesgue Space one, see
\cite{Fiorenza-Formica-Gogatishvili-DEA2018,fioformicarakodie2017,Kozachenko}.
 One can put
$$
\sigma_k := ||\eta||_k, \ k=1,2,\ldots,
$$
if, of course, $ \ \sigma_k \in (0,\infty)$.
 }
\end{remark}

\begin{corollary}
 Suppose that there exists a finite constant $ \ C < \infty \ $ and a
 non - negative r.v. $ \ \nu \ $ for which

$$
\sup_{k=0,1,\ldots} {\bf P} \left( \ \frac{|\eta_k|}{\sigma_k} > y \
\right) \le C \ {\bf P} (\nu > y), \ \ y \ge 1.
$$
Then the condition \eqref{T condition} is quite equivalent to the
following one

\begin{equation} \label{log moment}
{\bf E} \ln(e + \nu) < \infty.
\end{equation}

\end{corollary}
\noindent (See \cite{Rohatgi}).

\vspace{4mm}

 \ Note that the condition \eqref{log moment} is very weak; it is satisfied if, for
 example,

\begin{equation} \label{very weak}
\exists  \ b = {\const} > 0 \ : \  {\bf P} (\nu > z) \le C \ z^{-b},
\ \ z \ge 3.
\end{equation}

\vspace{4mm}

 \ Let us ground an opposite estimate for the radius of convergence.  We find quite analogously that if

\begin{equation} \label{q argument}
\forall q \in (0,1) \ \Rightarrow \sum_{k=1}^{\infty} {\bf P} \left( \ \frac{|\eta_k|}{\sigma_k} < q^k \ \right) < \infty,
\end{equation}
then

\begin{equation} \label{Non Gauss left radii}
\overline{\lim}_{n \to \infty} \sqrt[n]{|\eta_n|} \ge  \overline{\lim}_{n \to \infty} \sqrt[n]{\sigma_n}.
\end{equation}

\vspace{4mm}

 \ To summarize: \par

\vspace{5mm}

\begin{theorem}\label{theorem radius deterministic}
 Suppose that both the relations  \eqref{more general} and
\eqref{q argument} are satisfied. Then the radius $ \ R \ $ of
convergence for the power (Taylor) series \eqref{defin non Gauss
series} is deterministic almost surely and may be calculated by the
formula

\begin{equation} \label{Radii Taylor}
R = \frac{1}{\overline{\lim}_{n \to \infty} \sqrt[n]{\sigma_n}}.
\end{equation}

\vspace{4mm}
\end{theorem}

\begin{corollary}
 Suppose that there exists a non-negative r.v. $\, \tau \, $ for which $ \ \tau \ge 4 \ $
and, for some finite constant $ \ C < \infty \ $,  holds
$$
\sup_{k=0,1,\ldots} {\bf P} \left( \ \frac{\sigma_k}{|\eta_k|} > y \
\right) \le C \ {\bf P} (\tau > y), \ \ y \ge 1.
$$
Then the condition \eqref{q argument} is quite equivalent to the
following one

\begin{equation} \label{log moment at origin}
{\bf E} \ln(e + \tau) < \infty.
\end{equation}
\end{corollary}

\vspace{4mm}

\begin{example}
{\rm Evidently, if the r. v.-s $ \ \{\eta_k\} \ $ are independent,
then the radius of convergence of the series \eqref{defin non Gauss
series} is non-random.\par Let now  the sequence of coefficients $ \
\{\eta_k\} \ $ in  the \eqref{defin non Gauss  series}
 be such that

$$
{\bf P}(\forall k = 0,1,\ldots \ \Rightarrow \eta_k = 1) =
\frac{1}{2},
$$

$$
{\bf P}(\forall k = 0,1,\ldots \ \Rightarrow \eta_k = 2^{-k}) =
\frac{1}{2}.
$$

 \ On the other words, the infinite-dimensional vector of coefficients $ \ \eta = \{\  \eta_k \}, \ k = 0,1,2,\ldots
 $,
has the following distribution

$$
{\bf P}(\eta = (1,1,\ldots,1, \ldots) \ )= 1/2,
$$

$$
{\bf P}(\eta = (1,1/2, 1/4, \ldots, 1/2^k, \ldots) \ )= 1/2.
$$

 \ Then the radius of convergence $ \ R \ $  of the series \eqref{defin non Gauss  series} is a random variable
such that
$$
{\bf P} (R = 1) = {\bf P}(R = 2) = 1/2.
$$
}
\end{example}

 \vspace{3mm}

\begin{example}
{ \rm Let  $ \ \sigma_k = 1 \ $ and $ \ \tau \ $ be a positive
random variable such that $ \ \tau \ge e \ $ and

$$
{\bf P} (\tau \ge y) = \frac{1}{ \sqrt{\ln y}}, \ \ z \ge e.
$$
 \ Let us consider the following power series

$$
g_{\tau}(z) := \sum_{k=0}^{\infty} \tau \cdot z^k,
$$
i.e. here $ \ \eta_k = \tau, \ k = 0,1,\ldots.  \ $ The radius of
convergence is equal to one, despite that

$$
\forall Q > 1 \  \Rightarrow   \ \sum_{k=0}^{\infty} {\bf P} (\tau >
Q^k) = \infty.
$$
}
\end{example}

\vspace{5mm}

Let us ground the relations \eqref{rho f} and \eqref{beta f} of the
Theorem \ref{theorem order-type non random} in the general case,
i.e. non - Gaussian case.

 \vspace{4mm}

 \ Denote as before

$$
\rho \stackrel{def}{=} \overline{\lim}_{n \to \infty} \frac{n \ \ln n}{|\ln \sigma_n|}.
$$

\vspace{3mm}

 \ {\sc  Question:  under which conditions the following relation}

\vspace{3mm}

\begin{equation} \label{under non Gauss cond}
 \overline{\lim}_{n \to \infty} \frac{n \ \ln n}{|\ln \ |\eta_n| \ | \ } = \rho
\end{equation}

\vspace{3mm}

 \ {\sc  holds with probability one?} \par

\vspace{5mm}

 Namely one can assume, without loss of generality, $ \rho \in (0,\infty)$. Let $ \epsilon \in (0,\min(1,\rho/2))$ be an arbitrary constant.
 We have, for all the sufficiently
greatest values $ n \ge n_0 = n_0(\epsilon)$,

$$
\frac{n \ \ln n}{|\ln \sigma_n|} \le \rho(1 + \epsilon).
$$

 \ One can suppose $ \ \sigma_n \in (0,1); \ $ therefore

$$
0 < \sigma_n \le \frac{n^{-n \ \ln n}}{\rho^n (1 + \epsilon)^n} = n^{ \ -n/\rho(1 + \epsilon) \ }.
$$
 \ We deduce, from the direct definition of the variables $ \ \sigma_n \
 $,

\begin{equation*}
\begin{split}
{\bf P} \left( \ \frac{k \ \ln k}{|\ln |\eta_k| \ |} < \rho(1 + 2
\epsilon) \ \right) & = {\bf P}  \left( \frac{|\eta_k|}{\sigma_k} >
k^{\ k \epsilon/[\rho( 1 + \epsilon)(1 + 2 \epsilon) ] \ }
\right)\\
& \le T_k  \left( \ k^{\ k \epsilon/[\rho( 1 + \epsilon)(1 + 2
\epsilon) ] \ } \ \right), \ \ k \ge 1.
\end{split}
\end{equation*}

\vspace{4mm}

 \ We conclude again as a consequence, by virtue of lemma of Borel - Cantelli,
 the following result.

\vspace{3mm}

\begin{proposition}
If
\begin{equation} \label{sigma k}
\forall \delta > 0 \ \Rightarrow \sum_{k=1}^{\infty} T_k \left( \ k^{\delta k}   \ \right) < \infty,
\end{equation}
then, with probability one,
\begin{equation} \label{upp non Gauss}
\overline{\lim}_{n \to \infty} \left\{ \frac{n \ \ln n}{|\ln| \eta_n| \ |} \ \right\} \le \rho.
\end{equation}
\end{proposition}

\vspace{4mm}

\begin{corollary}
 Suppose that there exists a random variable $ \
\nu, \ \nu \ge 4,  \ $ and a constant $ \ C < \infty \ $ for which
 $$
 \sup_{k \ge 1} {\bf P} (|\eta_k| >z) \le C \ {\bf P}(\nu > z), \ \ z \ge 4.
 $$
The condition \eqref{upp non Gauss} is satisfied if, for instance,
$$
\forall \delta > 0 \ \Rightarrow \sum_{k=4}^{\infty} {\bf P} \left( \ \nu \ge k^{\delta k}  \ \right) < \infty.
$$
\end{corollary}
 \vspace{3mm}

 \ The last relation may be transformed  in turn as follows:

\begin{equation} \label{transf rel}
{\bf E} \left\{ \ \frac{\ln \nu}{\ln \ln \nu} \ \right\} < \infty.
\end{equation}

\vspace{4mm}

 \ Analogously can be justified the following fact. \par

 \vspace{5mm}

\begin{proposition}\label{prop lower bound} (Lower bound)

 Suppose that

 \begin{equation} \label{low non  Gauss}
\forall \delta > 0 \ \Rightarrow \sum_{k=1}^{\infty} {\bf P} \left( \  \frac{\sigma_k}{|\eta_k|} > k^{\delta \ k}  \ \right) < \infty.
 \end{equation}
Then, with probability one,
\begin{equation} \label{ result low non Gauss}
\overline{\lim}_{n \to \infty} \left\{ \frac{n \ \ln n}{|\ln|
\eta_n| \ |} \ \right\} \ge \rho.
\end{equation}

\end{proposition}

\vspace{4mm}

\begin{proof}
 \ Indeed, suppose for simplicity
$$
\frac{n \ln n}{|\ln \sigma_n|} = \rho = \const \in (0, \infty), \ \
n \ge 4,
$$
then
$$
\sigma_n = n^{ -n/\rho}.
$$
Let as above $ \ \epsilon \ $ be an arbitrary number from the (open)
interval $ \ (0, \min(1,\rho/2)). \ $ Consider a {\it sequence} of
events
$$
E_k(\epsilon) = \left\{\frac{k \ \ln k}{|\ln|\eta_k| \ |} \le \rho - \epsilon \   \right\}.
$$
We have
\begin{equation*}
{\bf P} (E_k(\epsilon)) \ = \ {\bf P} \left( \
\frac{|\eta_k|}{\sigma_k}  < k^{ - k \epsilon/(\rho(\rho - \epsilon)
) \ } \right) \ = \  {\bf P} \left( \ \frac{|\sigma_k|}{\eta_k}  >
k^{ k \epsilon/(\rho(\rho - \epsilon) ) \ } \right).
\end{equation*}

 It follows immediately, from the condition \eqref{low non  Gauss}, that

$$
\forall \epsilon \in (1,\min(1,\rho/2)) \ \Rightarrow \ \sum_{k=4}^{\infty} {\bf P} (E_k(\epsilon)) < \infty.
$$

 \ It remains to mention  once again the lemma of Borel - Cantelli. \par

\end{proof}
\vspace{4mm}

 \begin{remark}
{\rm  As above, the condition  \eqref{low non  Gauss} is satisfied
in turn if there exists a random variable $ \ \zeta, \ \zeta \ge 4,
\ $ for which

$$
\sup_{k \ge 4} {\bf P} \left(\ \frac{\sigma_k}{|\eta_k|} \ < z \ \right) \le C \cdot {\bf P} (\zeta > z), \ z \ge 4,
$$
and such that

$$
{\bf E} \left\{ \ \frac{\ln \zeta}{\ln \ln \zeta} \ \right\} < \infty.
$$
}
\end{remark}

\vspace{4mm}

 \begin{remark}
{\rm  Of course, if all the conditions of the Propositions 3.1 and
3.2 are satisfied, then with probability one

\begin{equation} \label{finally}
\overline{\lim}_{n \to \infty} \left\{ \frac{n \ \ln n}{|\ln| \eta_n| \ |} \ \right\} = \rho.
\end{equation}
}
\end{remark}

\vspace{5mm}

 \begin{remark}
{\rm   The mentioned problem was considered in many works, e.g.
\cite{Arnold1,Arnold2,Roters,Shapovalovska}, etc. As a rule, it was
investigated only the case when the coefficients $ \ \{\eta_k\} \ $
are independent or at last pairwise independent.}
\end{remark}

\vspace{5mm}

\section{Concluding remarks.}

\vspace{5mm}

 \ \hspace{4mm} {\bf A.}  \ Note that the classical methods offered by  V.M. Zolotarev  in \cite{Zolotarev} and, after by W. H\"oeffding in
 \cite{Hoffding}, bring more complicated computations. \par

\vspace{3mm}

 \ {\bf B.} The asymptotical behavior of the coefficients and, as a consequence, the asymptotical behavior of the function itself, allows to
 calculate the asymptotical distribution of its zeros, see e.g. the classical monograph on B.Ya. Levin (\cite{Levin}); see also
\cite{Hough,Peres,Sodin1,Sodin2}.

\vspace{6mm}

\vspace{0.5cm} \emph{Acknowledgement.} {\footnotesize The first
author has been partially supported by the Gruppo Nazionale per
l'Analisi Matematica, la Probabilit\`a e le loro Applicazioni
(GNAMPA) of the Istituto Nazionale di Alta Matematica (INdAM) and by
Universit\`a degli Studi di Napoli Parthenope through the project
\lq\lq sostegno alla Ricerca individuale\rq\rq}.\par

\end{document}